\documentclass[12pt,a4paper,twoside]{article}
\usepackage{a4}

\usepackage{amsmath, amssymb, amsthm}
\usepackage[dvips]{graphicx}
\def\<{\langle}
\def\>{\rangle}
\newcommand{\eps}{\varepsilon}
\def\dt{\partial_t\/}
\def\ZZ{\mathbb Z}
\def\RR{\mathbb R}

\def\calf{\mathcal{F}}
\def\call{\mathcal{L}}
\def\tr{\operatorname{Tr}}
\def\id{\operatorname{id}}
\def\Div{\operatorname{div}}
\newcommand\Ric{\operatorname{Ric}}
\newcommand\vol{\operatorname{vol}}

\newcommand{\eq}{\hspace*{-2mm}&=&\hspace*{-2mm}}

\newtheorem{cor}{Corollary}
\newtheorem{df}{Definition}
\newtheorem{ex}{Example}
\newtheorem{rem}{Remark}
\newtheorem{lem}{Lemma}
\newtheorem{prop}{Proposition}
\newtheorem{thm}{Theorem}
\newcounter{rembango}

\title{Extrinsic geometric flows\\ on foliated manifolds, II
\footnote{The work was supported by grant P-IEF, No. 219696 of Marie-Curie action.}}

\author{Vladimir Rovenski$^{(1,2)}$\thanks{E-mail: rovenski@math.haifa.ac.il
                                           }
       \ and
       Pawe\l\ Walczak$^{(2,3)}$\thanks{E-mail:
                                           pawelwal@math.uni.lodz.pl} \\
       \normalsize (1) Haifa University (2) Uniwersytet \L\'odzki
       (3) Instytut Matematyczny PAN}

\date{}

\begin{document}

\maketitle

\begin{abstract}
Extrinsic Geometric Flow (EGF) for a codi\-mension-one foliation
has been recently introduced by authors
as deformations of Riemannian metrics subject to quantities expressed in
terms of its second fundamental form.
In the paper we introduce soliton solutions to EGF and study their geometry
for totally umbilical foliations, foliations on surfaces,
and when the EGF is produced by the extrinsic Ricci tensor.

\vskip1mm
{\bf AMS Subject Classification:} 53C12; Secondary: 53C44

{\bf Keywords and Phrases:} geometric flow,  soliton, totally umbilical, principal curvatures, extrinsic Ricci curvature
\end{abstract}

\section*{Introduction}

The Hamilton's \textit{Ricci Flow} and the \textit{Mean Curvature Flow}
are among the central themes in recent Riemannian geometry.
The special \textit{soliton} solutions of these flows motivate the general analysis of the singularity formation.
The recent years have seen a growing interest in \textit{Geometric Flows} $\dt g_t = h(g_t)$ of different types.
The \textit{Extrinsic Geometric Flow (EGF)} has been recently introduced by authors \cite{rw0} as leaf-wise deformations of Riemannian metrics on a manifold $M$ equipped with a codi\-mension-one foliation $\calf$ subject to quantities expressed in terms of the second fundamental form of leaves.

 In the paper we introduce soliton solutions to EGF and study their geometry
for totally umbilical foliations, foliations on surfaces,
and when the EGF is produced by the extrinsic Ricci tensor.
Throughout the work,
{\bf
$(M^{n+1}, g)$ is a Riemannian space
with a codimension one transversely oriented foliation $\calf$,
$\nabla$ the Levi-Civita connection of~$g$,
$N$ the positively oriented unit normal to $\calf$,
$A:X\in T\calf\mapsto-\nabla_X N$ the Weingarten opera\-tor of the leaves,
which we extend to a $(1,1)$-tensor field on $TM$ by $A(N)=0$}.

\vskip1mm
The~definition $S(X_1, X_2)=S(\hat X_1,\hat X_2)\ (X_i\in TM)$ of the $\calf$-{\it truncated} $(0,2)$-{\it tensor field} $S$ ($\ \hat{}\,$ denotes the $T\calf$-component) will be helpful throughout the paper.
 Let $\hat g$ and $g^\perp$ be components of $g$ along $T\calf$ and $T\calf^\perp$, respectively.

Let $b:T\calf\times T\calf\to\RR$ be the second fundamental form of (the leaves of)
$\calf$ with respect to  $N$, and $\hat{b}$ its extension to
the $\calf$-truncated symmetric $(0, 2)$-tensor field on $M$.
In other words, $\hat{b}(N, \cdot)=0$ and $\hat{b}$ is dual to the extended Weingarten operator $A$. Denote by $\hat b_j$ the symmetric $(0,2)$-tensor fields on $M$
dual to powers $A^j$ of extended Weingarten operator,
\begin{equation*}
 \hat{b}_0(X, Y) =\hat g(X,Y),
 \quad
 \hat{b}_j(X, Y) = \hat g(A^j X, Y)\quad (j>0,\quad X,Y\in TM).
\end{equation*}
The symmetric functions $\tau_j=\tr A^j=\sum_{i=1}^n k_i^j\ (j\ge0)$, are called \textit{power sums} of the principal curvatures $\{k_i\}$ of $\calf$.
They can be expressed using \textit{elementary symmetric functions}
$\sigma_j=\sum\nolimits_{i_1<\ldots<i_j}k_{i_1}\cdot\ldots\cdot k_{i_j}\
(0\le j\le n),$ called \textit{mean curvatures}. Certainly, $\tau_{n+i}$ are not independent: using Newton formulae
\begin{eqnarray*}
 &&\hskip-12mm \tau_{j}-\tau_{j-1}\/\sigma_1+ \ldots
 +(-1)^{j-1}\tau_{1}\sigma_{j-1}+(-1)^j j\,\sigma_{j}=0 \quad (1\le j\le n),\\
 &&\hskip-12mm\tau_{j}-\tau_{j-1}\/\sigma_1+ \ldots+(-1)^n\,\tau_{j-n}\sigma_n=0
 \quad (j>n),
\end{eqnarray*}
one can express the functions $\tau_{n+i}\ (i>0)$ as polynomials of $\vec\tau:=(\tau_1, \ldots, \tau_n)$.

\begin{df}[see \cite{rw0}]\rm
Given functions $f_j\in C^{2}(\RR^{n})\ (0\le j<n)$ (at least one of them is not identically zero),
a family $g_t$ of Riemannian metrics on $(M,\calf, N)$ is called an \textit{Extrinsic Geometric Flow} (EGF), whenever
\begin{equation}\label{eq1}
 \dt g_t = h(b_t), \quad\text{where} \quad
 h(b_t)=\sum\nolimits_{j=0}^{n-1} f_j(\vec\tau)\,\hat{b}_j^{\,t}
\end{equation}
and $\hat b_j^{\,t}$ are $(0,2)$-tensors dual to powers $A^j_t$ of extended Weingarten operator $A_t$.
\end{df}

Indeed, EGF preserves $N$ unit and perpendicular to $\calf$, the $\calf$-compo\-nent of the vector does not depend on $t$.
EGF preserves the following classes of foliations:
 \textit{totally umbilical} ($A=\lambda\,\hat\id$),
 \textit{totally geodesic} ($A=0$) and \textit{Riemannian}
 ($\nabla_N N{=}\,0$).

Although (\ref{eq1}) consists of first order non-linear PDE's,
the corresponding power sums $\tau_i\ (i>0)$ satisfy the quasilinear system of PDE's
\begin{equation}\label{E-tauAk}
 \dt\tau_i+\frac{i}2\big\{\tau_{i-1}N(f_0)+\sum\nolimits_{j=1}^{n-1}\big[\,\frac{j f_j}{i+j-1} N(\tau_{i+j-1}) +\tau_{i+j-1}N(f_j)\big]\hskip-1pt\big\} = 0.
\end{equation}
In \cite{rw0}, we proved the local existence and uniqueness theorems for EGF and estimated the existence time of solutions for some cases.

For generic setting of $f_j$'s\ (while $f_0(0)=0$), the \textit{fixed points} of EGF
are totally geodesic ($A_t\equiv 0$) foliations only.
Several classes of foliations appear as fixed points of the flow $\dt g_t=f(\vec\tau)\,\hat g_t$ for special choices of $f$, for example:

-- foliations of constant $\tau_i$, when $f=\tau_i-c$ (minimal for $i=1$ and $c=0$),

-- totally umbilical foliations, when
   $f=n\,\tau _2-\tau _1^2 = \sum _{\,i<j}(k_i-k_j)^2$,

-- totally geodesic foliations, when $f(0)=0$, \ etc.

\vskip.5mm
Let $\mathrm{Diff}(M)$ be the diffeomorphism group of~$M$.
 Introduce the notation

-- $\mathcal{D}(\calf)$ the subgroup of $\mathrm{Diff}(M)$ preserving $\calf$;

-- $\mathcal{D}(\calf,N)$ the subgroup of $\mathrm{Diff}(\calf)$ preserving both, $\calf$ and $N$.

Recall that Fibration Theorem of D.\,Tischler (see for example \cite{cc})
states that the property of a closed manifold

\vskip.5mm
(i) \textit{$M$ admits a codimension one $C^1$-foliation invariant by a transverse flow};

\noindent
is equivalent to any of the following conditions:

\vskip.5mm
(ii) \textit{$M$ fibers over the circle $S^1$};

\vskip.5mm
(iii) \textit{$M$ supports a closed $1$-form of a class $C^1$ without singularities}.

\vskip.5mm\noindent
Taking into account the Theorem of R.\,Sacksteder, see for example \cite{cc},
one obtains in a class $C^2$ the following equivalent to any of (i),\,(ii) or (iii) condition:

\vskip.5mm
(iv) \textit{$M$ admits a codimension one foliation without holonomy}.

\begin{df}\label{D-ssEGsol}\rm
We say that a solution $g_t=\hat g_t\oplus g_t^\perp$ to (\ref{eq1})
is a \textit{self-similar Extrinsic Geometric Soliton} (EGS) on $(M,\calf,N)$ if there exist a smooth function
$\sigma(t)>0\ (\sigma(0)=1)$,
and a family of diffeomorphisms $\phi_t\in\mathcal{D}(\calf,N)$, $\phi_0=\id_M$, such that
\begin{equation}\label{E-solit1}
 \hat g_t=\sigma(t)\,\phi^*_t\,\hat g_0.
\end{equation}
\end{df}

By $\mathcal{M}(M,\calf,N)$, we denote the space of smooth Riemannian metrics on $M$ of finite volume with $N$ -- a unit normal to $\calf$.
Let $\RR_+$ acts on $\mathcal{M}$ by scalings along $T\calf$.
(By Remark~\ref{R-iii} in what follows, the Weingarten operator $A$ and the principal curvatures of $\calf$ are invariant under uniform scaling of the metric on~$\calf$).

The EGF may be regarded as a dynamical system on the quotient space $\mathcal{M}(M,\calf,N)/(\mathcal{D}(\calf,N)\times\RR_+)$.
The solutions to (\ref{eq1}) of the form (\ref{E-solit1})
correspond to \textit{fixed points} of the above dynamical system.

\vskip.5mm
\noindent
\textbf{Question}.
\textit{Given $(M,\calf,N)$, $N$ being a transversal vector field to $\calf$, ${\rm codim}\,\calf=1$, and $f_j\in C^{2}(\RR^{n})\ (0\le j<n)$, do there exist complete EGS metrics on $M$? If~they do exist, study their properties, classify them, etc.}

\vskip.5mm
The above Question is closely related to the basic problem (first posed by H.\,Gluck for geodesic foliations) in the theory of foliations.

\vskip.5mm
\noindent
\textbf{Problem} (see \cite{cw}).
\textit{Given (P), a property of Riemannian submanifolds
expressed in terms of the second fundamental form and its invariants, and
a foliated manifold $(M,\calf)$, decide if there exist (P)\,-metrics on
$(M,\calf)$, that is, Riemannian metrics such that all the leaves enjoy
the property~(P). If they do exist, study their properties, classify them
etc}.

\vskip.5mm
Some solutions to Problem for the Reeb foliation $\calf_{R}$ on a 3-sphere are the following:

(i)~\textit{there exist metrics for which $\calf_{R}$ is totally umbilic} \cite{bg},
 and

(ii)~\textit{there is no metric making $\calf_{R}$ minimal (or totally geodesic)}.

\vskip.5mm
In the paper we introduce and study EGS's as generalized fixed points of EGF's (see Theorems~\ref{T-0ERS}, \ref{T-1ERS}), and characterize them for some particular cases:
 totally umbilical foliations (Theorems~\ref{T-7umb}, \ref{T-Ntotumb}), foliations on surfaces (Corollaries~\ref{C-surf1},~\ref{C-surf2} and  Theorem~\ref{P-surf3}), for the EGF produced by extrinsic Ricci tensor (Theorem~\ref{T-ERS}).

\section{Preliminaries}
\label{sec:prelim}

It is well known that the Ricci tensor $\Ric$ is preserved by diffeomorphisms
$\phi\in{\rm Diff}(M)$ in the sense that $\phi^*(\Ric(g))=\Ric(\phi^*g)$,
hence $\Ric$ is an \textit{intrisic geometric tensor}.
In the same way, the second fundamental form $b$ of a foliation is invariant
under diffeomorphisms preserving the foliation,
hence $b$ is an \textit{extrinsic geometric tensor} for the class of foliations.
More precisely, we have the following

\begin{prop}\label{P-extprop}
Let $(M_i,\calf_i, g_i)\ (i=1,2)$ be Riemannian manifolds with codimension 1 foliations,
and $\phi:M_1\to M_2$ a diffeomorphism such that $\calf_1=\phi^{-1}(\calf_2)$ and
$g_1=\sigma\cdot\phi^*\hat g_2\oplus\phi^*(g_2^\perp)$ (for some positive $\sigma\in\RR$).
Then the Weingarten operators $A_i$ and the second fundamental forms $b^{(i)}$ satisfy
\begin{equation}\label{E-extinv1}
  b^{\,(1)}=\sigma\cdot\phi^* b^{\,(2)},\qquad
  A_1=\phi_*^{-1} A_2\,\phi_*, \qquad
  \vec{\tau}^{\,(1)}=\vec{\tau}^{\,(2)}\circ\phi,
\end{equation}
where $\vec\tau^{\,(i)}$ is the set of $\tau$-s for $A_i$.
For the tensor $h(b)=\sum\nolimits_{j=0}^{n-1} f_j(\vec\tau)\,\hat b_j$,
see Definition~\ref{D-ssEGsol}, we have
\begin{equation}\label{E-extinv-h}
  h(b^{\,(1)})=\sigma\cdot\phi^*(h(b^{\,(2)})).
\end{equation}
\end{prop}

\textbf{Proof} follows from the direct computation (see also Lemma~\ref{L-btAt} in what follows).
We include it here just for the convenience of a reader.

Notice that for any vector $X$ tangent to $\calf_1$, the vector $\phi_*X=d\phi(X)$ is tangent to $\calf_2$.
By definition of the metric $g_1$, we have $g_1(X,Y)=\sigma\cdot g_2(\phi_*X,\phi_*Y)$, where $X,Y\in T\calf_1$.
Let $N_2$ be a unit normal of $\calf_2$. Then $N_1=\phi_*^{-1}N_2$ is the unit normal of $\calf_1$.
In fact,
\begin{eqnarray*}
 g_1(X,N_1)\eq g_2(\phi_*X, N_2)=0\qquad (X\in T\calf_1),\\
 g_1(N_1,N_1)\eq\phi^*(g_2)(\phi_*^{-1}N_2,\phi_*^{-1}N_2)
 =g_2(\phi_*\phi_*^{-1}N_2,\phi_*\phi_*^{-1}N_2)=g_2(N_2,N_2).
\end{eqnarray*}
For the second fundamental form $b^{\,(2)}:T\calf_2\times T\calf_2\to\RR$ of $\calf_2$,
we obtain a symmetric $(0,2)$-tensor field $\phi^*b^{\,(2)}:T\calf_1\times T\calf_1\to\RR$.
Denote $\overset{i}\nabla$ the Levi-Civita connection of $(M_i,g_i)$.
For any vector fields $X,Y\subset T\calf_1$ we have
\begin{eqnarray*}
 (\phi^*b^{\,(2)})(X,Y)\eq b^{\,(2)}(\phi_*X, \phi_*Y)
 =-g_2(\overset{2}\nabla_{\phi_*X}N_2, \phi_*Y)
 =-g_2(\overset{2}\nabla_{\phi_*X}\phi_*N_1, \phi_*Y)\\
 \eq-g_2(\phi_*\overset{1}\nabla_{X}N_1, \phi_*Y)
 =-\phi^*(g_2)(\overset{1}\nabla_{X}N_1, Y) =\sigma^{-1}\cdot b^{\,(1)}(X,Y).
\end{eqnarray*}
Hence, the second fundamental form of a foliation is an \textit{extrinsic geometric tensor}.
Consequently, from the above we have
\[
 \phi_*(A_1 X)=A_2(\phi_* X)\quad (X\in T\calf_1) \ \ \Rightarrow \ \
 A_1=\phi_* A_2 \phi_*^{-1}.
\]
Taking the trace, we get the identity for $\tau$'s.
 For $h(b^{\,(i)})$ we have
\[
 \sum\nolimits_{j=0}^{n-1} f_j(\vec\tau^{\,(1)})\,\hat b^{\,(1)}_j(X, Y)=
 \sigma\sum\nolimits_{j=0}^{n-1} f_j(\vec\tau^{\,(2)}\circ\phi)\,\hat b^{\,(2)}_j(\phi_*X, \phi_*Y).
\]
From the above it follows (\ref{E-extinv-h}).
\qed

\vskip1.5mm
Next lemma is standard and we omit its proof.

\begin{lem}
\label{L-btAt}
Given a Riemannian manifold $(M,\,g=\hat g\oplus g^\perp)$ with a codimen\-sion-1 foliation $\calf$ and a unit normal $N$,
define a metric $\bar g = e^{\,2\/\varphi}\hat g\oplus g^\perp$ on $M$, where $\varphi:M\to\RR$ is a smooth function. Then
\begin{equation*}
 \bar b = e^{\,2\/\phi}(b -N(\varphi)\,\hat g),\qquad
 \bar A = A-N(\varphi)\,\widehat{\id}.
\end{equation*}
 \end{lem}

\begin{rem}\label{R-iii}\rm
If $\varphi$ is constant, then by Lemma~\ref{L-btAt} we have

\vskip.5mm
 (i)~ $\bar b=e^{\,2\/\varphi}\,b$ (the second fundamental forms);

\vskip.5mm
 (ii)~$\bar A = A$ (the Weingarten operators);

\vskip.5mm
 (iii)~ $\bar\tau_j=\tau_j$ (the power sums).
\end{rem}

Vector fields infinitesimally represent diffeomorphisms.
Let $\mathcal{X}(M)$ be Lie algebra of all vector fields on $M$ with the bracket operation. Introduce the notation

\vskip.5mm
-- $\mathcal{X}(\calf)$ the set of vector fields on $M$ preserving $\calf$,

\vskip.5mm
-- $\mathcal{X}(\calf,N)$ the set of vector fields on $M$ preserving both, $\calf$ and $N$.

\vskip.5mm\noindent
By the Jacobi identity, $\mathcal{X}(\calf)$ and $\mathcal{X}(\calf,N)$ are subalgebras of the Lie algebra $\mathcal{X}(M)$.
Moreover, for any $X\in\mathcal{X}(\calf)$ (or $X\in\mathcal{X}(\calf,N)$) there exists a family $\phi_t\in\mathcal{D}(\calf)$ (resp., $\phi_t\in\mathcal{D}(\calf,N)$)
such that $X={d\phi_t}/{dt}$ at $t=0$.

\begin{rem}\rm
(a) The following conditions are equivalent, \cite{wa}:
\begin{equation*}
 X\in\mathcal{X}(\calf) \ \Longleftrightarrow \ [X, T\calf]\subset T\calf.
\end{equation*}
If $\phi_t\in\mathcal{D}(\calf,N)$, then $\varphi_{t*} N = N\circ\varphi_t$.
For $X\in\mathcal{X}(\calf,N)$ generated by $\phi_t$, the above yields $\call_X N = 0$.
Putting these facts together, we conclude that
\begin{equation}\label{E-condbasicFN}
 X\in\mathcal{X}(\calf,N) \ \Longleftrightarrow \
 [X, T\calf]\subset T\calf\ \ {\rm and}\ \ [X, N] = 0.
\end{equation}

(b) By the above, a normal vector field $X=e^f N$ preserves $\calf$ if and only if
\begin{equation*}
 \nabla^\calf f= -\nabla_N N.
\end{equation*}
Here $\nabla^\calf\!f$ is the $\calf$-gradient of a function $f\in C^1(M)$.
Indeed, using (\ref{E-condbasicFN}), we get
\begin{equation*}
  0 = g([X, Y], N) = -e^f g(\,\nabla_N N, Y) + Y(e^f),\qquad \forall\ Y\perp N.
\end{equation*}
In particular, $N\in\mathcal{X}(\calf,N)$ if and only if $\calf$ is a Riemannian foliation.
\end{rem}

The \textit{Lie derivative} of a $(0,p)$-tensor $S$ w. r. to a vector field $X$ is given~by
\begin{equation}\label{E-LieDer}
 (\call_X S)(Y_1,\ldots, Y_p)=X(S(Y_1,\ldots, Y_p))
 -\sum\nolimits_{i=1}^p S(Y_1,\ldots,\call_X Y_i,\ldots,Y_p).
\end{equation}

\begin{lem}\label{L-LieXg}
For any vector fields $X\in\mathcal{X}(\calf,N)$ and $Y_i\in \Gamma(T\calf)$, we have
\begin{eqnarray*}
 (\call_X\,\hat g)(Y_1,Y_2)\eq \hat g(\nabla_{Y_1} X, Y_2) +\hat g(\nabla_{Y_2} X, Y_1),\\
 (\call_X\,\hat g)(N,Y_i)\eq\hat b_1(X, Y_i)-\hat g(\nabla_{X^\perp} N, Y_i),\qquad
 (\call_X\,\hat g)(N,N)=0.
\end{eqnarray*}
\end{lem}

\textbf{Proof}. Using $\nabla g=0$, we obtain
\begin{equation*}
 (\nabla_X\,\hat g)(N, N)=(\nabla_X \,\hat g)(Y_1, Y_2)=0,\quad
 (\nabla_X\,\hat g)(N, Y_i)=\hat b_1(X, Y_i) +\hat g(\nabla_{X^\perp} N, Y_i).
\end{equation*}
By above and (\ref{E-LieDer}), we have
\begin{eqnarray*}
 (\call_X \hat g)(Y_1,Y_2)\eq D_X(\hat g(Y_1,Y_2))-\hat g([X,Y_1], Y_2)-\hat g(Y_1, [X,Y_2])\\
 \eq\hat g(\nabla_{Y_1} X, Y_2) +\hat g(\nabla_{Y_2} X, Y_1),\\
 (\call_X\,\hat g)(N,Y_i) \eq-\hat g([X, N], Y_i)=b(X,Y_i)-\hat g(\nabla_{X^\perp} N, Y_i).
\end{eqnarray*}
 Similarly, $(\call_X\,\hat g)(N, N)=0$.
 Notice that $(\call_N\,\hat g)(N,Y_i)=-\hat g(\nabla_{N} N, Y_i)$.
\qed

\section{Introducing extrinsic geometric solitons}

Using Proposition~\ref{P-extprop}, we obtain

\begin{prop}\label{P-solitontau}
Let $g_t$ be a self-similar EGS (see Definition~\ref{D-ssEGsol}). Then
 \begin{eqnarray}\label{E-h-extend0}
  b_t\eq\sigma(t)\,\phi^*_t\,b_0,\qquad
  A_t=\phi_{\,t*}^{-1} A_0\,\phi_{\,t*},\qquad
  \vec\tau^{\,t}=\vec\tau^{\,0}\circ\phi_t,\\
\label{E-h-extend}
 h(b_t)\eq\sigma(t)\,\phi^*_t\,h(b_0).
\end{eqnarray}
\end{prop}

We are looking at what initial conditions give rise to self-similar EGS's.
Differentiating (\ref{E-solit1}) yields
\begin{equation}\label{E-solit1c}
 h(b_t)=\sigma'(t)\,\phi^*_t\hat g_0+\sigma(t)\phi^*_t(\call_{X(t)}\hat g_0),
\end{equation}
where $X(t)\in\mathcal{X}(\calf,N)$ is a time-dependent vector field generated by the family $\phi_t$.
  Since $h(b_t)=\sigma(t)\phi^*_t h(b_0)$, one may omit the pull-back in (\ref{E-solit1c}),
\begin{equation}\label{E-solit1d}
  h(b_0) =\sigma'(t)/\sigma(t)\,\hat g_0 +\call_{X(t)}\hat g_0.
\end{equation}

Motivated by the above, we accept the following

\begin{df}\label{D-EGsol}\rm
A pair $(g,X)$ consisting of a metric $g=\hat g\oplus g^\perp$ on $(M,\calf,N)$,
and a complete vector field $X\in\mathcal{X}(\calf,N)$ satisfying
for some $\eps\in\RR$ the condition
\begin{equation}\label{E-solit2}
 h(b) = \eps\,\hat g +\call_X\hat g\qquad
 {\rm where}\quad h(b)=\sum\nolimits_{0\le j<n}f_j(\vec\tau)\,\hat b_j
\end{equation}
is called a \textit{EGS structure}.
(We~say that $X$ is \textit{the vector field the EGS is flowing along}).
If~$X=\nabla f$ for some $f\in C^1(M)$, we have a \textit{gradient EGS structure}.
In~this case, $\frac12\call_{\nabla f}\,\hat g=\widehat{\rm Hess}_{g}f$
(the $\calf$-truncated hessian), and
\begin{equation}\label{E-solit3}
 h(b) =\eps\,\hat g + 2\,\widehat{\rm Hess}_{g} f.
\end{equation}
If~$X$ in (\ref{E-solit2}) is orthogonal to $N$, we have an \textit{$\calf$-EGS structure},
and if~$X$ is parallel to $N$, we have an \textit{$N$-EGS structure}.
For a {gradient $\calf$-EGS} one has $N(f)=0$.
For a {gradient $N$-EGS}, the function $f$ is constant along the leaves.
\end{df}

Next we observe that Definitions \ref{D-ssEGsol} and \ref{D-EGsol} are, in fact, equivalent.

\begin{thm}\label{T-0ERS}
(a) If $g_t$ is a self-similar EGS on $(M,\calf,N)$, then there exists a vector field $X\in\mathcal{X}(\calf,N)$
such that the metric $g_0$ solves~(\ref{E-solit2}).
(b)~Conversely, given vector field $X\in\mathcal{X}(\calf,N)$ and a~solution $g_0$ to (\ref{E-solit2}),
there is a function $\sigma(t)>0$ and a family of diffeomorphisms $\phi_t\in \mathcal{D}(\calf,N)$ such that a family of metrics $g_t$, defined by (\ref{E-solit1})
on $(M,\calf)$, is a solution to (\ref{eq1}).
\end{thm}

\textbf{Proof}.
(a) Recall that $\sigma(0)=1$ and $\phi_0=\id$.
Let $X=\frac d{dt}\phi_{t\,|t=0}\in\mathcal{X}(\calf,N)$ be the vector field generated by diffeomorphisms $\phi_t$.
Then we have
\begin{equation*}
 h(b_0)=\dt g_{t|\,t=0}=\dt\hat g_{t|\,t=0} =\sigma'(0)\hat g_0+\call_{X}\hat g_0.
\end{equation*}
This implies that $g_0$ and $X$ satisfy (\ref{E-solit2}) with $\eps=\sigma'(0)$.

(b) Conversely, suppose that $g_0$ satisfies (\ref{E-solit2}). Denote $\sigma(t)=e^{\eps\,t}$, hence $\sigma'(0)=\eps$.
Let $\psi_t\in\mathcal{D}(\calf,N)$ with $\psi_0=\id_M$ be a family of diffeomorphisms generated by $X$.
A~smooth family $g_t$ of $\calf$-truncated metrics on $M$, defined by $\hat g_t=\sigma(t)\psi_t^*\hat g_0$,
is of the form (\ref{E-solit1}).
 Moreover,
\begin{equation*}
 \dt g_t=\sigma'(t)\psi_t^*(\hat g_0)+\sigma(t)\psi_t^*(\call_{X}\hat g_0)
 =\frac{\sigma'(t)}{\eps}\psi_t^*(\eps\hat g_0+\call_{X}\hat g_0)=\frac{\sigma'(t)}{\eps}\,\psi_t^*(h(b_0)).
\end{equation*}
By Proposition~\ref{P-solitontau}, we have $f_j(\vec\tau^{\,0}\circ\psi_t)=f_j(\vec\tau^{\,0})$ for $t\ge0$,
and
\begin{eqnarray*}
 \psi_t^*(h(b_0))\eq\psi_t^*\big(\sum\nolimits_{j=0}^{n-1} f_j(\vec\tau^{\,0})\hat b_j^0\big)
 =\sum\nolimits_{j=0}^{n-1} f_j(\vec\tau^{\,0}\circ\psi_t)\,\psi_t^*\hat b_j^0\\
 \eq\sum\nolimits_{j=0}^{n-1} f_j(\vec\tau^{\,t})\,\sigma^{-1}(t)\,\hat b_j^t = \sigma^{-1}(t)\,h(b_t).
\end{eqnarray*}
Hence $h(b_t)=\sigma(t)\,\psi_t^*(h(b_0))$, and we conclude that
$\dt g_t=\frac{\sigma'(t)}{\sigma(t)\,\eps}\,h(b_t)=h(b_t)$.
\qed

\begin{rem}\rm
 If $(X,g)$ is a EGS structure, then  from
 Theorem~\ref{T-0ERS} and Proposition~\ref{P-solitontau}  it follows the constancy of $\tau$'s along~$X$:
\begin{equation}\label{E-Xtau}
 \nabla\!_X\tau_i = 0,\qquad 1\le i\le n.
\end{equation}
\end{rem}

\begin{ex}\rm
We will illustrate Theorem~\ref{T-0ERS}, case (b).
Let $(g_0,\,X)$ with $X=0$ be a EGS structure on $(M, \calf)$.
Then  $h(b_0)=\eps\hat g_0$ for some $\eps\in\RR$.
The family of leaf-wise conformal metrics $g_t=e^{\eps\,t}\hat g_0\oplus g^\perp_0$,
obviously satisfies PDE $\dt g_t=\eps\hat g_t$.
Using $\hat b_j^{\,t}=e^{\eps\,t}\,\hat b_j^{\,0}$ and $\vec\tau^{\,t}=\vec\tau^{\,0}$
(see Remark~\ref{R-iii}), we obtain $h(b_t)=e^{\eps\,t}h(b_0)=e^{\eps\,t}\,\eps\hat g_0=\eps\hat g_t$.
Hence $\dt g_t=h(b_t)$, and  $g_t$ is a self-similar EGS.
\end{ex}

\begin{thm}[Canonical form]\label{T-1ERS}
Let a self-similar EGS $g_t$ be unique among soliton solutions to (\ref{eq1}) with initial metric $g_0$.
Then there is a 1-parameter family of diffeomorphisms $\psi_t\in\mathcal{D}(\calf,N)$ and
a constant $\eps\in\{-1,0,1\}$ such that
\begin{equation}\label{E-solit2c}
 \hat g_t = (1+\eps\,t)\,\psi_t^*\hat g_0.
\end{equation}
\end{thm}

\noindent
 The cases $\eps=-1,0,1$ in (\ref{E-solit2c}) correspond to \textit{shrinking, steady}, or \textit{expanding}~EGS.

\vskip1mm
\textbf{Proof} is similar to the one of \cite[Proposition 1.3]{ck2}.
For convenience of a reader we provide it here in the case $\sigma''(0)\ne0$.
From (\ref{E-solit1d}) it follows
\begin{equation}\label{E-solit1nice}
   h(b_0)= (\log\sigma)'(t)\,\hat g_0+\call_{X(t)}\hat g_0,
 \end{equation}
where $X(t)\in\mathcal{X}(\calf,N)$ is a family of vector fields such that
$X(t)=d\phi_{t}/dt$. Differentiating (\ref{E-solit1nice}) with respect to $t$ gives
\begin{equation}
 (\log\sigma)''(t)\,\hat g_0 +\call_{{X}'(t)}\hat g_0=0.
\end{equation}
Let $Y_0=-{X}'(0)/(\log\sigma)''(0)$. We then have $\call_{Y_0}\hat g_0=\hat g_0$.
Substituting this into (\ref{E-solit1d}), we have for all $t$
\[
 h(b_0)=\call_{(\log\sigma)'(t)Y_0+X(t)}\,\hat g_0.
\]
Put $X_0=(\log\sigma)'(0)Y_0+X(0)$. Then $h(b_0)=\call_{X_0}\hat g_0$.
Let $\psi_t\in\mathcal{D}(\calf,N)$ be a family of diffeomorphisms generated by $X_0$.
We will check that
\begin{equation}\label{E-prop3}
 \tilde g_t=\psi^*_t\hat g_0\oplus(\psi^*_t g_0)^\perp.
\end{equation}
is the EGF with the same initial conditions $g_0$, and that it is a steady soliton
(i.e., $\sigma(t)=1$ for all $t$).
Indeed, differentiating (\ref{E-prop3}), we have by (\ref{E-h-extend}),
\[
 \dt\tilde g_t=\psi^*_t(\call_{X_0}\hat g_0)=\psi^*_t(h(b_0))=h(\psi^*_t b_0)=h(\tilde b_t).
\]
Thus $g_t=\tilde g_t$, by uniqueness assumption for EGS solutions to
our flow with initial metric $g_0$.
By replacing $\phi_t$ by $\psi_t$ we have $\sigma(t)\equiv1$ in (\ref{E-solit1}).
\qed

\begin{rem}\rm
Equation (\ref{E-solit2}) yields a rather strong condition on the EGS structure $(g, X)$.
For example, contracting (\ref{E-solit2}) with $g$ (tracing) and using the identity
$\tr\call_X \hat g=2\Div_\calf X$ (see Lemma~\ref{L-LieXg}) yields
\begin{equation}\label{E-solit4a}
 \tr_g h(b) = n\,\eps + 2\Div_\calf X.
\end{equation}
For a gradient EGS, (\ref{E-solit4a}) means
 $\tr_g h(b) = n\,\eps + 2\,\Delta_\calf f$.
\end{rem}

\begin{prop}
(a) Let $(g,X)$ be a EGS structure on $(M,\calf)$. If a leaf $L$ is compact, then
\begin{equation}\label{E-intepsL}
 n\,\eps = {\int_L \tr_g h(b)\,d\vol_{g,L}}/{\vol(L,g)}.
\end{equation}
(b) Let $M$ be compact and either $X\perp N$ and $\nabla_N N=0$ or $X||N$ and $\tau_1=0$.~Then
\begin{eqnarray}\label{E-eps}
 n\,\eps = {\int_M \tr_g h(b)\,d\vol_g}/{\vol(M,g)}.
\end{eqnarray}
In particular, if $\calf$ is Riemannian with minimal leaves,
then (\ref{E-eps}) holds.
\end{prop}

\textbf{Proof}.
Integrating (\ref{E-solit4a}) and using the Divergence Theorem we obtain (\ref{E-intepsL}).
 From $\Div(f N)=f\Div N +N(f)=-\tau_1 f+N(f)$,
by the Divergence theorem, $\int_{M}\Div(f N)\,d\vol=0$,
we obtain
\begin{equation*}
 \int_M N(f)\,d\vol=\int_M \tau_1 f\,d\vol.
\end{equation*}
Next, we have $g(\nabla_N X, N)=N g(X, N)- g(X,\nabla_N N)$, and by the above,
\begin{eqnarray*}
 \int_M g(\nabla_N X, N)\,d\vol_g =
 \int_M \big( \tau_1\,g(X, N)- g(X, \nabla_N N)\big) d\vol_g.
\end{eqnarray*}
Therefore, in case (b), integrating of (\ref{E-solit4a}) implies (\ref{E-eps}).
 If $\calf$ is Riemannian with minimal leaves, then $\tau_1=0$ and $\nabla_N N=0$.
\qed

\begin{prop}\label{P-hb}
Equation (\ref{E-solit2}) for EGS with $X=\mu\,N$ $(\mu:M\to\RR_+)$ reads as
\begin{equation*}
 h(b) = \eps\hat g - 2\,\mu\,\hat b_1.
\end{equation*}
For Riemannian foliations we, certainly, have $\mu\equiv 1$.
\end{prop}

\textbf{Proof}.
From (\ref{E-solit2}) and Lemma~\ref{L-LieXg} (for $X=\mu\,N$) we obtain
\begin{equation*}
 (\call_{\mu N}\,\hat g)(Y_1,Y_2) = -2\,\mu\,\hat b_1(Y_1,Y_2),\qquad Y_i\perp N.
 \qquad\qed
\end{equation*}

\begin{ex}\rm
Let $h(b)=\hat b_1$ (i.e., $f_j=\delta_{j1}$). Any metric making $(M,\calf)$ a Riemannian foliation, is a steady EGS with $X=N$ (unit normal).
Indeed, $N\in\mathcal{X}(\calf,N)$ for a bundle-like metric $g$ (see Lemma~\ref{L-LieXg}),
and we have $h(b)=\frac12\,\call_N\hat g$.
\end{ex}

\section{Totally umbilical EGS structures}
\label{sec:tumb}

Recall that the EGF preserves total umbilicity of $\calf$.

\begin{prop}\mbox{\rm\cite{rw0}}\label{P-08}
 Let $g_t\ (0\le t<\eps)$ be the EGF (\ref{eq1}) on $(M,\calf,N)$.
 If $\calf$ is totally umbilical for $g_0$, then $\calf$ is totally umbilical for any $g_t$.
\end{prop}

Let $\calf$ be a totally umbilical foliation on $(M ,g)$ with the normal curvature $\lambda$ (not identically zero).
We have $A=\lambda\,\hat\id$ and $\tau_j=n\lambda^j$. Hence EGF is given by
\begin{equation}\label{E-rotsymflow}
 \dt g_t=\psi(\lambda_t)\,\hat g_t,
\end{equation}
where
\begin{equation}\label{E-psi}
  \psi(\lambda)=\sum\nolimits_{j=0}^{n-1}f_j(n\lambda,n\lambda^2,\ldots,n\lambda^n)\,\lambda^j.
\end{equation}
In this case, $\lambda_t$ obeys the quasilinear PDE, see (\ref{E-tauAk}),
\begin{equation}\label{E-EGF-lambda}
 \dt\lambda_t+\frac12 N(\psi(\lambda_t))=0,
\end{equation}
and the EGS structure equations (\ref{E-solit2}) reduce to the PDE
\begin{equation}\label{E-psieps-2}
 \psi(\lambda) = \eps + (2/n)\,\Div_\calf\tilde X.
\end{equation}
If $(g, X)$ is a EGS structure with totally umbilical metric, by (\ref{E-Xtau}) we have
\begin{equation}\label{E-Xlambda}
  X(\lambda) = 0.
\end{equation}
 The EGF on a surface $(M^2,g)$, foliated by curves $\calf$, is given by (\ref{E-rotsymflow}),
 where $\lambda=\tau_1$ is the geodesic curvature of the curves-leaves,
 and $\psi=f_0\in C^2(\RR)$.
 The EGS equations (\ref{E-solit2}) on $(M^2,\calf)$ reduce to the PDE, see (\ref{E-psieps-2}) with $n=1$,
\begin{equation}\label{E-psieps}
 \psi(\lambda) = \eps + 2\,\Div_\calf X.
\end{equation}

\begin{rem}\rm
Let $\calf$ be a codimension-1 totally umbilical foliation with the metric
\begin{equation}\label{E-bfconf1}
 g_t = e^{2 f_t}\hat g_0\oplus g_0^\perp,
\end{equation}
where $f_t:M\to\RR\ (f_0=0)$ are smooth functions. We claim that
\begin{equation}\label{eq11-B}
 2\,\dt f_t =\psi\big(\lambda - N(f_t)\big).
\end{equation}
Indeed, by Lemma~\ref{L-btAt}, $A_t=A_0 - N(f_t)\,\widehat{\id}$, hence, $\lambda_t=\lambda - N(f_t)$.
By Lemma~\ref{L-btAt}, we also have
 $b_t=e^{2\,f_t}(\lambda-N(f_t))\,\hat g_{0} =(\lambda-N(f_t))\,\hat g_{t}$.
Similarly, $b_t^{\,j}=(\lambda-N(f_t))^j\,\hat g_{t}$.
Differentiating (\ref{E-bfconf1}) yields
$
 \psi(\lambda_t)\hat g_t=h(b_t)=\dt g_t = (2\,\dt f_t)\,\hat g_t.
$
Hence, $2\,\dt f_t=\psi(\lambda_t)$ that gives us the non-linear PDE (\ref{eq11-B}).

One can solve explicitly for $f_t$ only the particular case $f=c_1\lambda+c_2$ of the problem
$\dt g_t=\psi(\lambda_t)\,\hat g_t$, when (\ref{eq11-B}) becomes linear of the form
\begin{equation*}
 2\,\dt f_t + N(f_t) = c_1\,\lambda_t+c_2.
\end{equation*}
In general, the non-linear PDE (\ref{eq11-B}) is difficult to solve, and we apply EGF approach, see \cite{rw0}:
first we find $\lambda_t$ from (\ref{E-EGF-lambda}), then find $g_t$ from (\ref{E-rotsymflow}).
\end{rem}

\begin{ex}\label{Ex-solit1}\rm
$\phantom{.}$

(a) Let $\calf$ be a~totally umbilical foliation  on $(M, g_0)$ with $\lambda={\rm const}$
(if $M$ is compact then $\lambda=0$ by known integral formula $\int_M\lambda\,f\vol=0$, see, for example, \cite{rw1}).
Using $h(b_0)=\psi(\lambda)\hat g_0$, see (\ref{E-rotsymflow}), we conclude that $g_0$ with $X=0$ and $\eps=\psi(\lambda)$ is a EGS structure.
Moreover, $g_t=e^{\psi(\lambda)\,t}\,\hat g_0\oplus g^\perp_0$ is a self-similar EGS with $\phi_t=\id_{M}$ and $\sigma(t)=e^{\psi(\lambda)\,t}$.

For a specific case of a totally geodesic foliation (i.e., $\lambda\equiv0$, if such $g_0$ exists on $(M,\calf)$),
$g_0$ with $X=0$ and $\eps=\psi(0)$ is a EGS structure,
and the family of totally geodesic metrics $g_t=e^{\psi(0)\,t}\,\hat g_0\oplus g^\perp_0$ is also a self-similar EGS.

\vskip1mm
(b) Let $(g,X)$ be a EGS structure on $(M,\calf)$.
If $\calf$ is totally umbilical (with the normal curvature $\lambda$), then $X$ is a leaf-wise conformally Killing field: $\call_X\hat g=(\psi(\lambda)-\eps)\,\hat g$.
If~$\calf$ is totally geodesic (hence $\psi=f_0({\bf 0})$), then $X$ is an \textit{infinitesimal homothety along leaves},
$\call_X\hat g=C\,\hat g$, with the factor $C=f_0({\bf 0})-\eps$.
In particular, $X$ is a \textit{leaf-wise Killing field} when $\eps=f_0({\bf 0})$.
This happens, for example, when $M$ is a surface of revolution in $M^{n+1}(c)$ foliated by parallels, see Example~\ref{Ex-surfrev}.

(c) Consider biregular foliated coordinates $(x_0,x_1)$ on a surface $M^2$ (see \cite[Section~5.1]{cc}).
Since the coordinate vectors $\partial_0,\partial_1$ are directed along $N$ and $\calf$, respectively,
the metric has the form $g=g_{00}\,dx_0^2+g_{11}\,dx_1^2$.
Recall that $h(b)=\psi(\lambda)\,g_{11}$.
By \cite[Lemma~4]{rw0} with $n=1$, we have $\lambda=-\frac1{2\sqrt{g_{00}}}(\log g_{11})_{,0}$.
Let~$X=X^0\partial_0+X^1\partial_1\in\mathcal{X}(\calf,N)$. Using $g_{01}=0$, we obtain
\begin{eqnarray*}
 \Div_\calf X=g(\nabla_{\partial_1} X, \partial_1)
 = (\,\partial_1(X^1) +X^0\Gamma^1_{01} +X^1\Gamma^1_{11})\,g_{11},
\end{eqnarray*}
where
 $\Gamma^1_{01}=\frac12(\log g_{11})_{,0}$ and $\Gamma^1_{11}=\frac12(\log g_{11})_{,1}$.
Hence, (\ref{E-psieps}) has the form
\begin{equation*}
 \psi(\lambda) -\eps = 2\,\partial_1(X^1)g_{11} +X^0 g_{11,0}+X^1 g_{11,1}.
\end{equation*}
From the condition $[X, \partial_1]\perp\partial_0$, see (\ref{E-condbasicFN}), and
 $[X,\ \partial_1] =-\partial_1(X^0)\partial_0-\partial_1(X^1)\partial_1$
we obtain $\partial_1(X^0)=0$.
Next, from the condition $[X, N]=0$, see (\ref{E-condbasicFN}), and
\begin{equation*}
 [X,\ g_{00}^{-\frac12}\partial_0\,] = (X(g_{00}^{-\frac12}) -g_{00}^{-\frac12}\,\partial_0(X^0))\partial_0 +g_{00}^{-\frac12}\,\partial_0(X^1)\partial_1
\end{equation*}
we obtain
\[
 \partial_0(X^1)=0,\qquad
 \partial_0(X^0)= -\frac12\,X(\log g_{00}).
\]
\end{ex}

\begin{df}[see \cite{nir}]\rm
Denote the torus $\RR^{n+1}/\ZZ^{n+1}$ by ${\rm T}^{n+1}$. For $v\in\RR^{n+1}$, let $R^t_v(x):=x+tv$ be the flow on ${\rm T}^{n+1}$ induced by a ''constant" vector field $X_v$.
We say $v\in \RR^{n+1}$ is \textit{Diophantine}
if there is $s>0$ such that $\inf\limits_{u\in\ZZ^{n+1}\setminus\{0\}} |\<u,v\>|\cdot\|u\|^s>0$, where
$\<\ ,\ \>$ and $\|\cdot\|$ are the Euclidean inner product and the norm in $\RR^{n+1}$.
  When~$v$ is Diophantine, we call $R_v$ a \textit{Diophantine linear flow}.
\end{df}

\begin{thm}\label{T-7umb}
Let $\calf$ be a codimension-one totally umbilical foliation (with a unit normal $N$)
on a torus $({\rm T}^{n+1}, g),\ n>0$.
Suppose that $X$  is a smooth unit vector field on $T\calf$ with the properties
\begin{equation}\label{E-Xcond12}
 (i)~\nabla_XX\,\perp\,T\calf,\qquad
 (ii)~R(X,Y)Y\in T\calf\qquad(Y\in T\calf).
\end{equation}
If $X$-flow is conjugate (by a homeomorphism) to a Diophantine liner flow $R_v$, then
for any function $\psi$ of a class $C^2$, see (\ref{E-rotsymflow}), (\ref{E-psi}),
there exists $f:{\rm T}^{n+1}\to\RR$ such that $(g, f X)$ is a EGS structure, see (\ref{E-psieps-2}).
\end{thm}

\textbf{Proof}.
For a totally umbilical foliation with the normal curvature
$\lambda$, the Weingarten operator is conformal: $A=\lambda\,\hat\id$. By (\ref{E-Xcond12})$_{ii}$ and the Codazzi equation
$$(\nabla_X A)Y-(\nabla_Y A)X=R(X,Y)N^\top,$$
we have $\lambda=const$ along $X$-curves.
By (\ref{E-Xcond12})$_{i}$, the $X$-curves are $\calf$-geodesics.

Let $(g,\,\tilde X)$ be a EGS structure with $\tilde X=fX$.
From (\ref{E-psieps-2}) and the known identity $\Div_\calf fX = f\Div_\calf X+ X(f)$
it follows that $\Div_\calf \tilde X = X(f)$.
We are looking for solution of PDE, see (\ref{E-psieps}),
\begin{equation*}
 \psi(\lambda)-\eps = (2/n)\,X(f).
\end{equation*}
Since $X$-flow is conjugate to a Diophantine linear flow, by Kolmogorov theorem,
see \cite{nir}, the above PDE has a solution
$(f,\eps)\in C^\infty({\rm T}^{n+1})\times\RR$.
\qed

\vskip1.5mm
From Theorem~\ref{T-7umb} it follows

\begin{cor}\label{C-surf1}
Let a unit vector field $X$ on a torus $({\rm T}^2, g)$ defines a foliation $\calf$ by curves of constant geodesic curvature $\lambda$. If $X$-flow is conjugated to a Diophantine liner flow $R_v$,
then for any function $\psi$ of a class $C^2$, see (\ref{E-rotsymflow}), (\ref{E-psi}),
there exists $f:{\rm T}^2\to\RR$ such that $(g, f X)$ is a EGS structure, see (\ref{E-psieps}).
\end{cor}

Notice that if $\psi\in C^2(\RR)$ then the following function belongs to $C^1$:
\begin{equation}\label{E-mulambda}
 \mu= \left\{
 \begin{array}{c}
 -\frac n2(\psi(\lambda)-\psi(0))/\lambda,\quad \lambda\ne0,\\
 \ -\frac n2\psi\,'(0),\qquad\qquad\qquad \lambda=0.\\
 \end{array}\right.
\end{equation}

\begin{thm}\label{T-Ntotumb}
Let $\calf$ be a totally umbilical foliation on $(M, g)$ with the normal curvature $\lambda$,
and the function $\psi\in C^2(\RR)$ given in (\ref{E-psi}) satisfies $\psi'\ne0$.
Then the following properties are equivalent:

\vskip1mm
(1) the normal curvature of $\calf$ satisfies $N(\lambda)=0$;

\vskip1mm
(2) $(g, \mu N)$, for some function $\mu$, is a EGS structure, see (\ref{E-psieps-2}),

(indeed, one may take $\mu$ as in (\ref{E-mulambda}) and $\eps=\psi(0)$).
\end{thm}

\textbf{Proof}.
$(1)\Rightarrow(2)$:
The EGS equations (for a totally umbilical metric $g$ and the vector field $X=\mu\,N$) are, see Proposition~\ref{P-hb} and  (\ref{E-Xlambda}),
\begin{equation}\label{E-psilambdaeps}
 \psi(\lambda) - \eps = -(2/n)\,\mu\,\lambda,\qquad X(\lambda)=0.
\end{equation}
For $\eps=\psi(0)$ and $\mu$ given in (\ref{E-mulambda}), the above (\ref{E-psilambdaeps}) are satisfied.
Hence, by Definition~\ref{D-EGsol}, $(g, \mu N)$ satisfy (\ref{E-psieps-2}).

$(2)\Rightarrow(1)$:
Using Definition~\ref{D-EGsol}, (\ref{E-Xlambda}) and $\psi'\ne0$,
we obtain the equality $\mu N(\lambda)=0$ with $\mu$ given in (\ref{E-mulambda}).
Consider an open set $\Omega=int\{p\in M: \mu=0\}$. Indeed, $N(\lambda)=0$ on $M\setminus\Omega$.
By~(\ref{E-mulambda}), we have $\psi(\lambda)=\psi(0)$ and hence $N(\psi(\lambda))=\psi'(\lambda)N(\lambda)=0$ on $\Omega$.
Since $\psi'\ne0$, we have $N(\lambda)=0$ on $\Omega$.
From the above we conclude that $N(\lambda)=0$ on $M$.
\qed

\vskip1mm
From Theorem~\ref{T-Ntotumb} it follows

\begin{cor}\label{C-surf2}
Let $\calf$ be a foliation (by curves) on a surface $(M^2, g)$
and $\psi\in C^2(\RR)$ is given in (\ref{E-psi}) and satisfies $\psi'\ne0$.
Then the following properties are equivalent:

\vskip1mm
(1) the geodesic curvature $\lambda$ of $\calf$ satisfies $N(\lambda)=0$;

\vskip1mm
(2) $(g, \mu N)$, for some function $\mu$, is a EGS structure, see (\ref{E-psieps}),

(one may take $\mu$ as in (\ref{E-mulambda}) with $n=1$ and $\eps=\psi(0)$).
\end{cor}

\begin{ex}[\it Non-Riemannian EGS on double-twisted products]\rm
 Let $M=M_1\times M_2$ be the product (with the metric  $\tilde g=g_1\oplus g_2$ and Levi-Civita connection $\tilde\nabla$) of a circle $M_1=S^1$ with the canonical metric $g_1$ and a compact Riemannian manifold $(M_2,g_2)$.
 Let $f_i:M\to \RR$ be a positive differentiable function,
 $\pi_i:M\to M_i$ the canonical projection,
 $\pi^\perp_i:TM\to\ker \pi_{3-i}$ the vector bundle projection, for $i=1,2$.
 The \textit{metric of a double-twisted product}
 $M_1\times_{(f_1\times f_2)} M_2$ is
\begin{equation*}
  g(X,Y) = f_1 g_1(\pi_{1*}X, \pi_{1*}Y) + f_2 g_2(\pi_{2*}X, \pi_{2*}Y),\quad X,Y\in TM,
\end{equation*}
i.e., $g=f_1 g_1\oplus f_2 g_2$.
The Levi-Civita connection $\nabla$ of $g$ obeys the relation \cite{pr}
\begin{equation*}
 \nabla_X Y = \tilde\nabla_X Y + \sum\nolimits_i \big(g(\pi_{i*}X, \pi_{i*}Y) U_i -g(X, U_i)\pi_{i*}Y -g(Y, U_i)\pi_{i*}X\big).
\end{equation*}
Both foliations $M_1\times\{p_2\}$ and $\{p_1\}\times M_2$ are totally umbilical with the mean curvature vectors $H_1=\pi^\perp_{2}U_1$
and $H_2=\pi^\perp_{1}U_2$, respectively, where $U_i = -\nabla(\log f_i)$.
 This property characterizes the double-twisted product, see \cite{pr}.

Suppose that dim\,$M_1=1$, then EGF preserves the above double-twisted product structure.
The mean curvature of the foliation $\calf:=\{p_1\}\times M_2$ is constant along
the fibers $M_1\times\{p_2\}$ (i.e., $N$-curves)
if and only if $\pi^\perp_{1}U_2$ is a function of~$M_2$.
In this case, due to Theorem~\ref{T-Ntotumb}, $\calf$ admits a EGS structure with $X\,||\,N$.
\end{ex}

\begin{thm}\label{P-surf3}
 Let $(g,X)$ be a EGS structure on a closed surface $M^2$ foliated  by curves $\calf$ of the geodesic curvature $\lambda$,
 and let $\psi\in C^2(\RR)$, see (\ref{E-psi}), satisfies $\psi'\ne0$.
 If $X\,||\,N$ then $X$-curves are closed and define a fibration $\pi:M^2\to S^1$,
 and $\calf$ is the suspension of a diffeomorphism $f:S^1\to S^1$.
 Moreover, if $\psi(\lambda)=-2\lambda+c$ holds then  $(g,X)$ is the EGS structure (with $\eps=c$) for any metric $g\in\mathcal{M}$ satisfying $N(\lambda)=0$,
 otherwise $\lambda=0$ (i.e., $\calf$ is a geodesic foliation).
\end{thm}

\textbf{Proof}.
Assume an opposite, then the foliation $\calf_N$ has a limit cycle. Since $M$ is compact, there is
a domain $\Omega\subset M^2$ bounded by closed $N$-curves (which are limit cycles).
By (\ref{E-Xlambda}), $\lambda=const$ along $N$-curves. Since there are limiting leaves, $\lambda=const$ on $\Omega$.
From the relation $\Div N=-\lambda$ and the Divergence Theorem
\begin{equation*}
 \int_{\Omega}\Div N\,d\,\vol=\int_{\partial\Omega}\<N, \,\nu\>\,d\,\omega,
\end{equation*}
where $\nu$ is the outer normal to the boundary $\partial\Omega$ (hence $\nu\perp N$), we conclude that $\lambda=0$ on $\Omega$, hence $\calf_{N}$ is a Riemannian foliation, -- a contradiction.

By the classification theorem for foliations on closed surfaces,  see \cite{g},
all the $X$-curves are closed and define a fibration $\pi:M^2\to S^1$.
By (\ref{E-psilambdaeps}) with $\mu=1$ we have the following.
If $\psi(\lambda)\ne-2\lambda+c$ then $\lambda=const$ on $M$.
Hence, using the integral formula $\int_M\lambda\,d\vol=0$, we conclude that
$\lambda=0$.
In this case, by Lemma~\ref{L-btAt}, any $\calf$-geodesic metric on $(M,\calf,N)$
has the form $\bar g=(\pi^{-1}\circ\sigma)\hat g\oplus g^\bot$, where $\sigma:S^1\to\RR$ is an arbitrary smooth function.
\qed

\section{Rotational symmetric metrics}
\label{sec:rotsymm}

Notice that the EGF preserves \textit{rotational symmetric metrics}
 \begin{equation*}
 g = dx_0^2+\varphi^2(x_0)\,ds_{n}^2\quad
 \mbox{where }\ ds_{n}^2\ \  \mbox{is a metric of curvature 1}.
\end{equation*}
 The $n$-parallels $\{x_0=c\}$ form a  Riemannian totally umbilical foliation $\calf$ with the unit normal field $N=\partial_0$. In this case,
 EGF equation has the form (\ref{E-rotsymflow}) discussed in Example~\ref{Ex-solit1}.
 Any leaf-wise Killing field $X\perp N=\partial_0$ provides the EGS structure on $M^{n+1}$ with the rotational symmetric metric $g$.

Assuming $\hat g_t=\varphi^2_t\,\hat g_0$ and using (\ref{E-rotsymflow}), we obtain
$\lambda_t= -(\varphi_t)_{,0}/\varphi^2_t$ \ and
\begin{eqnarray*}
 \dt\varphi_t=\frac12\,\psi(\lambda_t)\,\varphi_t
 \quad \Rightarrow \quad
 \varphi_t=\varphi_0\,\exp\big(\frac12\int_0^t\psi(\lambda_t)\,dt\big).
\end{eqnarray*}
In particular case of $\psi(\lambda)=\lambda$, we get the linear PDE
 $\dt\lambda+\frac12 N(\lambda)=0$
representing the ``unidirectional wave motion" along any $N$-curve $\gamma(s)$,
\begin{equation}\label{eq-unmotion}
 \lambda_t(s)=\lambda_0(s-t/2).
\end{equation}
If, in addition, $\lambda_0=C\in\RR$, then $\lambda_t=C$
and $\varphi_t=\varphi_0\,\exp(\frac12\,t\,\psi(C))$.

\begin{figure}[ht]
\begin{center}
\includegraphics[scale=.3,angle=0,clip=true,draft=false]{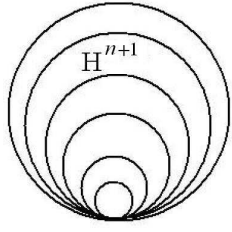}
\caption{\small Horosphere foliation.}
\label{F-002}
\end{center}
\end{figure}

Rotational symmetric metrics with $\lambda=\mbox{const}$ exist on hyperbolic space
${\mathbb H}^{n+1}$ with horosphere foliation.
On the Poincar\'{e} $(n+1)$-ball $B$ the leaves of such Riemannian totally umbilical foliations are Euclidean $n$-spheres tangent $\partial B$ (Fig.~\ref{F-002}).

The unit normal field $N$ and a leaf-wise Killing field $X\perp N$ provide $\calf$-EGS structures with $g$ as above.
Notice that there are no totally umbilical foliations with $\lambda=\mbox{const}\ne0$ on compact (or of finite volume, see \cite{rw1}) manifolds.

Some of rotational symmetric metrics come from hypersurfaces of revolution in space forms.
Evolving them by EGF yields deformations of hypersurfaces of revolution foliated by $n$-parallels.

\begin{ex}[EGS on hypersurfaces of revolution]\label{Ex-surfrev}\rm
Revolving the graph of $x_1=f(x_0)$ about the $x_0$-axis of $\RR^{n+1}$,
we get the hypersurface  $M^n:\ f^2(x_0)=\sum_{i=1}^n x_i^2$
foliated by $(n-1)$-spheres $\{x_0=c\}$ (parallels) with the induced metric
\begin{equation}\label{eq_g-rev}
 g=(1+{f'(x_0)}^{\,2})\,dx_0^2+f^2(x_0)\,\sum\nolimits_{i=1}^n dx_i^2.
\end{equation}

(i) Revolving a line $\gamma_0: x_1=\tan\beta\,x_0$ about the $x_0$-axis,
we build the cone $C_0: (\tan\beta\,x_0)^2=\sum_{i=1}^n x_i^2$,
with the metric $g_0=dx_0^2+(x_0\sin\beta)^2\,\sum_{i=1}^n dx_i^2$.
Hence $\varphi_0=x_0\sin\beta$ and $\lambda_0(x_0)=-2/x_0$.
Applying the EGF $\dt g_t=\lambda_t\hat g_t$, by (\ref{eq-unmotion}) we obtain
$\lambda_t(x_0)=-\frac2{x_0-t/2}$.
The rotational symmetric metric
$g_t=dx_0^2+(x_0-t/2)^2\sin^2\beta\,\sum_{i=1}^n dx_i^2$ appears on
the same cone translated across the $x_0$-axis,
$C_t: (x_0-t/2)^2\tan^2\beta=\sum_{i=1}^n x_i^2$.
Any leaf-wise Killing field $X\perp N$ provides the EGS structure on $M^n$ with the induced metric $g$.

(ii) Let us find a curve $y=f(x)>0$ such that the metric (\ref{eq_g-rev})
on the surface of revolution $M^n: \sum_{i=1}^n x_i^2=f^2(x_0)$ has $\lambda_0\,{=}\,\mbox{const}\,{=}\,1$. Using $\lambda_0=\frac1{f(x_0)}\sin\phi$, where $\tan\phi=f'(x_0)$, we get ODE
 $\,\frac{|f'|}{f\sqrt{1+(f')^2}}=1 \Rightarrow \frac{dx_1}{dx_0}=\frac{x_1}{\sqrt{4+x_1^2}}$.
\begin{figure}[ht]
\begin{center}
\includegraphics[scale=.4,angle=0,clip=true,draft=false]{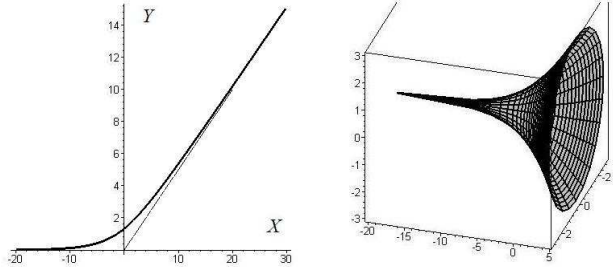}
\caption{\small a) Graph of $\gamma$.
\ \ b) Hypersurface of revolution (of $\gamma$) with $\lambda=const$.}
\label{F-001}
\end{center}
\end{figure}
The solution~is
 $\gamma:\ x_0=\log\frac{\sqrt{4+x_1^2}-2}{\sqrt{4+x_1^2}+2}+\sqrt{4+x_1^2}+C$,
 where $C\in\RR$.
The hypersurface $M^n\subset\RR^{n+1}$ looks like a pseudosphere (for $n=2$ see Fig.~\ref{F-001} and \cite[Example 4]{rw0}),
but for $x_0\to\infty$ it is asymptotic to the cone $(x_0+C)^2=\sum_{i=1}^n x_i^2$.
The~2-dimensional sectional curvature is $K(\partial_0,\partial_1)=-\frac{1}{(x_1^2+2)^2}<0$, and $\lim\limits_{x_1\to\pm\infty}K=0$.

As for horosphere foliation on ${\mathbb H}^{n}$,
the normal $N$ to parallels and any leaf-wise Killing field $X\perp N$
compose EGS structures on $(M^n,\,g)$.
\end{ex}

\section{Extrinsic Ricci flow}
\label{sec:extRic}

The {extrinsic Riemannian curvature tensor} $\mbox{Rm}^{\rm\,ex}$ of $\calf$
is, roughly speaking, the difference of the curvature tensors of $M$ and of the leaves.
More precisely, due to the Gauss formula, we have
$$
 \mbox{Rm}^{\rm\,ex}(Z, X)Y = g(A X,Y)A Z - g(A Z, Y)A X\quad {\rm for}\quad X,Y,Z\in T\calf.
$$
The \textit{extrinsic Ricci flow} is defined by
\begin{equation*}
 \dt g_t=-2\Ric^{\rm\,ex}_t,
\end{equation*}
where $\Ric^{\rm\,ex}(X,Y)=\tr\mbox{Rm}^{\rm\,ex}(\cdot,X)Y=\tau_1\hat b_1-\hat b_2$
is the \textit{extrinsic Ricci tensor}. For $n=2$, we have $\Ric^{\rm\,ex}=\sigma_2\,\hat g$.
Hence, $-2\Ric^{\rm\,ex}$ relates to $h(b)$ of (\ref{eq1}) with $f_1=-2\,\tau_1,\,f_2=2$ (others $f_j=0$).

\begin{ex}\rm
Foliations satisfying $\Ric^{\rm\,ex}=0$ are fixed points of extrinsic Ricci flow, they have the property $A(A-\tau_1\id)=0$, or,
  $k_j(k_j-\tau_1)=0\ (1\le j\le n)$ for the eigenvalues $k_j$ of $A$.
Since $\tau_1=\sum_j k_j$, from above it follows $k_j=0$ for all $j$.
Hence, extrinsic Ricci flat foliations are totally geodesic foliations only.
\end{ex}

In order to extend the set of solutions
we define the \textit{normalized EGF}~by
\begin{equation}\label{eq1-n}
  \dt g_t = h(b_t) -({\rho_t}/n)\,\hat g_t\quad{\rm with}\quad
  \rho_t={\int_M \tr A_h\,d\vol_t}/{\,{\rm Vol}(M, g_t)}.
\end{equation}

\begin{df}\rm
We call $g_t{\in}\,\mathcal{M}$ on a compact $M$ a \textit{normalized extrinsic Ricci flow}~if
\begin{equation}\label{eq-extRn}
  \dt g_t = -2\Ric_t^{\rm\,ex} +({\rho_t^{\rm\,ex}}/n)\,\hat g_t,\quad
  \rho_t^{\rm\,ex}=-{2\int_M \Ric_t(N,N)\,d\vol_t}/{\,{\rm Vol}(M,g_t)}.
\end{equation}
A metric $g$ on $(M, \calf)$ is \textit{extrinsic Einstein} if
 $\Ric^{\rm\,ex}={\rho}/({2\,n})\cdot\hat g$ for some $\rho\in\RR$.
Follow Definition~\ref{D-EGsol}, an \textit{extrinsic Ricci soliton structure}
is a pair $(g,X)$ of a metric $g$ on $(M,\calf)$, and a complete field $X\in\mathcal{X}(\calf,N)$
satisfying for some $\eps\in\RR$
\begin{equation*}
 -2\Ric^{\rm\,ex} = \eps\,\hat g +\call_X\hat g.
\end{equation*}
\end{df}

\begin{rem}\rm
 To explain $\rho_t^{\rm\,ex}$ in (\ref{eq-extRn}), we find
the \textit{extrinsic scalar curvature}:
$
 R^{\rm\,ex}=\tr\Ric^{\rm\,ex}=\tr(\tau_1\,A-A^2)=\tau_1^2-\tau_2=2\,\sigma_2.
$
By the integral formula $\int_M(2\,\sigma_2-\Ric(N,N))\,d\vol=0$ (see \cite{rw1}) we find
\begin{equation*}
 \int_M R^{\rm\,ex}\,d\vol=\int_M\Ric(N,N)\,d\vol.
\end{equation*}
Substituting this into (\ref{eq1-n}) instead of $\int_M \tr A_h\,d\vol$, we obtain $\rho_t^{\rm\,ex}$ of (\ref{eq-extRn}).
 Hence, extrinsic Einstein foliations are fixed points of the flow (\ref{eq-extRn}).
\end{rem}

Codimension one totally umbilical foliations with $\lambda=const$
are extrinsic Einstein foliations.
Clearly, an extrinsic Ricci soliton structure with $X=0$ contains extrinsic Einstein metric.

\begin{rem}\rm
A~codimension one foliation $(M,g)$ will be called \textit{CPC}
(\textit{constant principal curvatures}) if the principal curvatures of leaves are constant.

(a) (Non-)\,normalized EGF preserves CPC property of foliations, see \cite{rw0}.
Let such flow on $(M,\calf)$ starts from a CPC metric.
From $N(\tau_i)=0$, $N(f_j(\vec\tau))=0\ (j<n)$ and
(\ref{E-tauAk}) we conclude that $\tau_i$ do not depend on~$t$.

\vskip1mm
(b) Let $(G,g)$ be a compact Lie group with a left invariant metric $g$.
Suppose that the corresponding Lie algebra has a codimension one subspace $V$
such that $[V,V]\subset V$. Then $V$ determines a CPC foliation on $(G,g)$.
\end{rem}

\begin{thm}\label{T-ERS}
Let $(g,X)$ be an extrinsic Ricci soliton structure on $(M^n, \calf)\ (n>2)$,
and $X$ a leaf-wise conformal Killling field
(i.e., $\call_X \hat g = \mu\,\hat g$).

\hskip-3mm
(i) Then there are at most two distinct principal curvatures at any point $p\in M$.

\hskip-3mm
(ii) Moreover, if $\mu$ is constant along the leaves, then $\calf$ is CPC foliation.
\end{thm}

\textbf{Proof}. (i) Since $(\Ric^{\rm\,ex})^{\sharp}=-A(A-\tau_1\,\hat\id)$
and $\call_X \hat g = \mu\,\hat g$, we obtain the equality $A(A-\tau_1\id)=r\,\hat\id$ with
$r=\frac12(\eps+\mu)$, that yields equalities for the principal curvatures $k_j$,
\begin{equation*}
  k_j(k_j-\tau_1)=r\quad \forall j.
\end{equation*}
Hence each $k_j$ is a root of a quadratic polynomial $P_2(k)=k^2-\tau_1 k -r$.
The roots are real if and only if $\tau_1^2+4\,r\ge0$.
In the case $r>-\tau_1^2/4$ we have two distinct roots
$\bar k_{1,2}=(\tau_1\pm\sqrt{\tau_1^2+4\,r})/2$.
Let $n_1\in(0,n)$ eigenvalues of $A$ are equal to $\bar k_{1}$ and others to $\bar k_{2}$.
From $\tau_1=n_1\bar k_1+n_2\bar k_2$ and $n=n_1+n_2$ we obtain
\begin{equation}\label{eq-n1R}
  n_2-n_1={(n-2)\,\tau_1}/(\tau_1^2+4\,r)^{1/2}\ \in\mathbb{Z}.
\end{equation}
If $n_2=n_1$ then $\tau_1=0$ and $k_{1,2}=\pm\sqrt{r}$, otherwise,
$\tau_1^2=\frac{4\,r}{q^2-1}$ for $q=\frac{n-2}{n_2-n_1}$.

(ii) Assume that $\mu$ is constant along the leaves, then $r=const$ on any leaf.
In view of $n>2$, a continuous function $\tau_1: M\to\RR$ has values in a discrete set, hence it is constant.
Then all $k_j$'s (from both sets) are constant on $M$.
For $n$ even, (\ref{eq-n1R}) admits a particular solution:
 $\tau_1=0$, and $n/2$ principal curvatures $k_j$ equal to $\sqrt r$, others to $-\sqrt r$.
\qed

\vskip1mm
  By the above, extrinsic Einstein foliations satisfy the equality $A(A-\tau_1\id)=r\id$,
where $r=\rho/(2n)$. Hence, from Theorem~\ref{T-ERS} we obtain

\begin{cor}\label{C-EEinstf}
Let $\calf$ be a foliation with extrinsic Einstein metric $g$.

\hskip-3.5mm
(i) If $M$ is compact then $\calf$ is a fixed point of the normalized extrinsic Ricci~flow.

\hskip-3.5mm
(ii) If $n>2$ then $\calf$ is CPC foliation with $\le 2$ distinct principal curvatures.
\end{cor}


\begin{thebibliography}{WWW}

\bibitem[{BG}]{bg}
 M.\,Brunella, and E.\,Ghys, {\it Umbilical foliations and transversely holomorphic flows}.
 J. Differential Geom., 41(1) (1995) 1--19

\bibitem[{CC}]{cc}
 A.\,Candel and L.\,Conlon, \textit{Foliations, I}, AMS, Providence, 2000

\bibitem[{CW}]{cw} M.\,Czarnecki and P.\,Walczak,
Extrinsic geometry of foliations, 149--167, in ``\textit{Foliations 2005}", World Scientific Publ., NJ, 2006

\bibitem[{CK2}]{ck2} B.\,Chow et al,
\textit{The Ricci Flow: Techniques and Applications, Part I}, AMS, 2007

\bibitem[{G}]{g}
C.\,Godbillon, Dynamical Systems on Surfaces, Springer Verlag, Berlin -- Heidelberg -- New York, 1983

\bibitem[{Nir}]{nir}
 L.\,Nirenberg, \textit{Topics in Nonlinear Functional Analysis}. AMS, 2001

\bibitem[{PR}]{pr}
R.\,Ponge, and H.\,Reckziegel, \textit{Twisted products in pseudo-Riemannian geometry},
Geometriae Dedicata, 48 (1993) 15--25

\bibitem[{RW0}]{rw0}
V.\,Rovenski, and P.\,Walczak,
\textit{Extrinsic geometric flows on foliated manifolds, I}, (2010) 34 pp.
arXiv.org/math.DG/0001007\,v1 (2010)

\bibitem[{RW1}]{rw1}
V.\,Rovenski, and P.\,Walczak,
\textit{Integral formulae on foliated symmetric spaces},
Preprint, Univ. of Lodz, Fac. Math. Comp. Sci., 2007/13, 27\,pp. (2007)

\bibitem[{Wa}]{wa} P.\,Walczak, {\it Mean curvature invariant foliations}, Illinois J. Math., 37 (1993) 609--623

\end{thebibliography}
\end{document}